\documentclass[a4paper]{article}

\usepackage[tableposition=top,font=small,labelfont=bf,format=hang]{caption}
\usepackage{subfig}

\usepackage[latin1]{inputenc}
\usepackage[T1]{fontenc}
\usepackage{ae,aecompl}
\usepackage{amsmath,amsfonts}

\usepackage{booktabs}
\usepackage[sort&compress,numbers]{natbib}
\usepackage[dvips,colorlinks=false]{hyperref}
\hypersetup{
pdfauthor = {Lars Eirik Danielsen and Matthew G. Parker},
pdftitle = {On the Classification of All Self-Dual Additive Codes over GF(4) of Length up to 12},
pdfkeywords = {Self-dual codes, Graphs, Local complementation}
}
\usepackage[dvips]{graphicx}

\DeclareMathOperator{\GF}{GF}
\DeclareMathOperator{\wt}{wt}
\DeclareMathOperator{\Aut}{Aut}
\DeclareMathOperator{\tr}{Tr}

\usepackage{amsthm}
\newtheorem{thm}{Theorem}
\theoremstyle{definition}
\newtheorem{defn}[thm]{Definition}
\newtheorem{exmp}[thm]{Example}

\title{On the Classification of All Self-Dual Additive Codes over GF(4) of Length up to 12}

\author{Lars Eirik Danielsen\thanks{The Selmer Center, Department of Informatics, University of Bergen, 
PB 7800, \mbox{N-5020} Bergen, Norway.\hfill
\texttt{\{\href{mailto:larsed@ii.uib.no}{larsed},\href{mailto:matthew@ii.uib.no}{matthew}\}@ii.uib.no}\hfill
\texttt{http://www.ii.uib.no/\~{}\{\href{http://www.ii.uib.no/\~larsed}{larsed},\href{http://www.ii.uib.no/\~matthew}{matthew}\}}} \and Matthew G. Parker\footnotemark[1]}

\begin{document}

\date{February 17, 2006}
\maketitle

\begin{abstract}
We consider additive codes over $\GF(4)$ that are self-dual with respect to the
Hermitian trace inner product. Such codes have a well-known interpretation as quantum codes
and correspond to isotropic systems.
It has also been shown that these codes can be represented
as graphs, and that two codes are equivalent if and only if the corresponding graphs are
equivalent with respect to local complementation and graph isomorphism.
We use these facts to classify all codes of length up to 12, where previously only 
all codes of length up to 9 were known.
We also classify all extremal Type~II codes of length 14.
Finally, we find that the smallest Type~I and Type~II codes with trivial automorphism
group have length 9 and 12, respectively.\\[\baselineskip]
\emph{Keywords:} Self-dual codes; Graphs; Local complementation
\end{abstract}

\section{Introduction}

An \emph{additive} code, $\mathcal{C}$, over $\GF(4)$ of \emph{length}~$n$ is an
additive subgroup of $\GF(4)^n$. $\mathcal{C}$~contains $2^k$ codewords for some $0 \le k \le 2n$,
and can be defined by a $k \times n$ \emph{generator matrix}, with entries from $\GF(4)$,
whose rows span $\mathcal{C}$ additively. $\mathcal{C}$~is called an $(n,2^k)$ code.
We denote $\GF(4) = \{0,1,\omega,\omega^2\}$, where $\omega^2 = \omega + 1$.
\emph{Conjugation} of $x \in \GF(4)$ is defined by $\overline{x} = x^2$.
The \emph{trace map}, $\tr : \GF(4) \to \GF(2)$, is defined by
$\tr(x) = x + \overline{x}$.
The \emph{Hermitian trace inner product} of two vectors over $\GF(4)$ of length~$n$, 
$\boldsymbol{u} = (u_1,u_2,\ldots,u_n)$ and $\boldsymbol{v} = (v_1,v_2,\ldots,v_n)$, is given by
\begin{equation}
\boldsymbol{u} * \boldsymbol{v} = \tr(\boldsymbol{u} \cdot \overline{\boldsymbol{v}}) 
= \sum_{i=1}^n \tr(u_i \overline{v_i}) = \sum_{i=1}^n (u_i v_i^2 + u_i^2 v_i) \pmod{2}.
\end{equation}
Note that $\boldsymbol{u} * \boldsymbol{v}$ is also
the number (modulo 2) of places where $\boldsymbol{u}$ and $\boldsymbol{v}$ have different non-zero values.
We define the \emph{dual} of the code $\mathcal{C}$ with respect to
the Hermitian trace inner product, $\mathcal{C}^\perp = \{ \boldsymbol{u} \in \GF(4)^n \mid
\boldsymbol{u}*\boldsymbol{c}=0 \text{ for all } \boldsymbol{c} \in \mathcal{C} \}$.
$\mathcal{C}$~is \emph{self-orthogonal} if $\mathcal{C} \subseteq \mathcal{C}^\perp$.
It has been shown that self-orthogonal additive codes over $\GF(4)$
can be used to represent \emph{quantum error-correcting codes}~\cite{calderbank}.
If $\mathcal{C} = \mathcal{C}^\perp$, then $\mathcal{C}$ is \emph{self-dual}
and must be an $(n,2^n)$ code.
Self-dual additive codes over $\GF(4)$ correspond to zero-dimensional quantum codes,
which represent single quantum states.
If the code has high minimum distance, the corresponding quantum state is highly \emph{entangled}.

The \emph{Hamming weight} of $\boldsymbol{u}$, denoted $\wt(\boldsymbol{u})$,
is the number of nonzero components of $\boldsymbol{u}$.
The \emph{Hamming distance} between $\boldsymbol{u}$ and $\boldsymbol{v}$
is $\wt(\boldsymbol{u} - \boldsymbol{v})$.
The \emph{minimum distance} of the code $\mathcal{C}$ is the minimal Hamming distance
between any two distinct codewords of $\mathcal{C}$. Since $\mathcal{C}$ is an additive code,
the minimum distance is also given by the smallest nonzero weight of any codeword in $\mathcal{C}$.
A code with minimum distance~$d$ is called an $(n,2^k,d)$ code.
The \emph{weight distribution} of the code $\mathcal{C}$ is the sequence
$(A_0, A_1, \ldots, A_n)$, where $A_i$ is the number of codewords of weight~$i$.
The \emph{weight enumerator} of $\mathcal{C}$ is the polynomial
\begin{equation}
W(x,y) = \sum_{i=0}^n A_i x^{n-i} y^i
\end{equation}

We distinguish between two types of self-dual additive codes over $\GF(4)$. A code is 
of \emph{Type~II} if all codewords have even weight, otherwise it is of \emph{Type~I}.
It can be shown that a Type~II code must have even length.
Bounds on the minimum distance of self-dual codes were given by Rains and Sloane~\cite[Theorem~33]{rains}.
Let $d_I$ be the minimum distance of a Type~I code of length~$n$. Then $d_I$ is upper-bounded by
\begin{equation}
d_I \le \begin{cases}
2 \left\lfloor \frac{n}{6} \right\rfloor + 1, \quad \text{if } n \equiv 0 \text{ } (\text{mod } 6) \\ 
2 \left\lfloor \frac{n}{6} \right\rfloor + 3, \quad \text{if } n \equiv 5 \text{ } (\text{mod } 6) \\ 
2 \left\lfloor \frac{n}{6} \right\rfloor + 2, \quad \text{otherwise.}
\end{cases}
\end{equation}
There is a similar bound on $d_{II}$, the minimum distance of a Type~II code of length~$n$,
\begin{equation}
d_{II} \le 2 \left\lfloor \frac{n}{6} \right\rfloor + 2.
\end{equation}
A code that meets the appropriate bound is called \emph{extremal}. It can be shown
that extremal Type~II codes must have a unique weight enumerator.
Rains and Sloane~\cite{rains} also used a linear programming bound, and showed
that extremal codes do not exist for all lengths. For instance, there is no self-dual $(13,2^{13},6)$ code.
If a code has highest possible minimum distance, but is not extremal, it is called \emph{optimal}.
An interesting open problem is whether there exists a Type~II $(24,2^{24},10)$ code.

A \emph{linear} code, $\mathcal{C}$, over $\GF(4)$ which is self-dual with respect to the
\emph{Hermitian inner product}, i.e., $\boldsymbol{u} \cdot \overline{\boldsymbol{v}}=0 
\text{ for all } \boldsymbol{u},\boldsymbol{v} \in \mathcal{C}$,
is also a self-dual \emph{additive} code with respect to the \emph{Hermitian trace inner product}.
However, most of the self-dual additive codes are not linear.
Only Type~II codes can be linear, since self-dual linear codes over $\GF(4)$ must contain codewords
of even weight only. It follows that the set of Hermitian self-dual linear codes over $\GF(4)$ is a subset of 
the set of Type~II self-dual additive codes over $\GF(4)$.

\begin{exmp}\label{ex:hexacode}
The unique extremal $(6,2^6,4)$ code, also known as the \emph{Hexacode}, has a generator matrix
\[
\left(
\begin{array}{cccccc}
1&0&0&1&\omega&\omega\\
\omega&0&0&\omega&\omega^2&\omega^2\\
0&1&0&\omega&1&\omega\\
0&\omega&0&\omega^2&\omega&\omega^2\\
0&0&1&\omega&\omega&1\\
0&0&\omega&\omega^2&\omega^2&\omega
\end{array}
\right).
\]
This code has weight enumerator $W(x,y) = x^6 + 45 x^2 y^4 + 18 y^6$. It is therefore of Type~II, and it
can be verified that it is also a linear code.
\end{exmp}

Two self-dual additive codes over $\GF(4)$, $\mathcal{C}$ and $\mathcal{C}'$, are
\emph{equivalent} if and only if the codewords of $\mathcal{C}$ can be mapped onto the codewords 
of $\mathcal{C}'$ by a map that preserves self-duality.
Such a map must consist of a permutation of coordinates (columns of the generator matrix),
followed by multiplication of coordinates by nonzero elements from $\GF(4)$, followed by possible conjugation of 
coordinates. For a code of length~$n$, there is a total of $6^n n!$ such maps.
The 6 possible transformations given by scaling and conjugation of a coordinate 
are equivalent to the 6 permutations of the elements $\{1,\omega,\omega^2\}$ in the coordinate.
A map that maps $\mathcal{C}$ to $\mathcal{C}$ is called an \emph{automorphism} of $\mathcal{C}$.
All automorphisms of $\mathcal{C}$ make up an \emph{automorphism group},
denoted $\Aut(\mathcal{C})$. The number of distinct codes equivalent to $\mathcal{C}$ is then
given by $\frac{6^n n!}{|\Aut(\mathcal{C})|}$.
By summing the sizes of all equivalence classes, we find the total number of distinct codes
of length~$n$, denoted $T_n$.
It was shown by Höhn~\cite{hohn} that $T_n$ is also given by the \emph{mass formula},
\begin{equation}\label{eq:mass}
T_n = \prod_{i=1}^n (2^i+1) = \sum_{j=1}^{t_n} \frac{6^n n!}{|\Aut(\mathcal{C}_j)|},
\end{equation}
where the sum is over all equivalence classes.
Similarly, the total number of distinct Type~II codes of length~$n$ is given by
\begin{equation}\label{eq:mass2}
T_n^{\text{II}} = \prod_{i=0}^{n-1} (2^i+1) = \sum_{j=1}^{t_n^{\text{II}}} \frac{6^n n!}{|\Aut(\mathcal{C}_j)|},
\end{equation}
where the sum is over the equivalence classes of Type~II codes.
By assuming that $|\Aut(\mathcal{C}_j)| = 1$ for all $j$ in Eq.~\eqref{eq:mass}, we 
get a lower bound on $t_n$, the number of inequivalent codes of length $n$.
\begin{equation}\label{eq:bound}
t_n \ge \left\lceil\frac{\prod_{i=1}^n (2^i+1)}{6^n n!}\right\rceil
\end{equation}
A similar bound on $t_n^{\text{II}}$ can be derived from Eq.~\eqref{eq:mass2}.

We can use the computational algebra system \emph{Magma}~\cite{magma} to find the automorphism group 
of a code. Since, at this time, Magma has no explicit function for calculating the automorphism group of an additive code,
we use the following method, described by Calderbank~et~al.~\cite{calderbank}.
We map the $(n,2^k)$ additive code $\mathcal{C}$ over $\GF(4)$ to the $[3n,k]$ binary linear code
$\beta(\mathcal{C})$ by applying the map $0 \mapsto 000$,  $1 \mapsto 011$, $\omega \mapsto 101$,
$\omega^2 \mapsto 110$ to each generator of $\mathcal{C}$. 
We then use Magma to find $\Aut(\beta(\mathcal{C})) \cap \Aut(\beta(\GF(4)^n))$, which will be isomorphic
to $\Aut(\mathcal{C})$.

If we are given $t_n$ inequivalent codes of length $n$, i.e., one code from each equivalence class,
it is relatively easy to calculate the automorphism group size of each code, as described above,
and verify that the mass formula defined by Eq.~\eqref{eq:mass} gives the correct value.
But to actually find a set of $t_n$ inequivalent codes, or just the value of $t_n$, is a hard problem.
All self-dual additive codes over $\GF(4)$ of length~$n$ were first
classified, up to equivalence, by Calderbank~et~al.~\cite{calderbank} for $n \le 5$ and
by Höhn~\cite{hohn} for $n \le 7$. Höhn also classified all Type~II codes of length 8.
Using a different terminology, the codes of length~$n$ were implicitly classified by
Hein, Eisert, and Briegel~\cite{hein} for $n \le 7$ and by Glynn~et~al.~\cite{glynnbook} for $n \le 9$.
These classifications were not verified using the mass formula defined by
Eq.~\eqref{eq:mass}. Gaborit et~al.~\cite{gaborit,gaborit2} have classified all extremal codes
of length 8, 9, and 11, and all extremal Type~II codes of length 12. 
Bachoc and Gaborit~\cite{bachoc} classified all extremal Type~II codes of length 10,
and they also showed that there are at least 490 extremal Type~II codes of length 14
and gave a partial result on the unicity of the extremal Type~II code of length 18.
A review of the current status of the classification of various types of self-dual codes
is given by Huffman~\cite{huffman}.

In this paper, we will give a complete classification of all codes of length up to 12,
and all extremal Type~II codes of length 14.
But first, in Section~\ref{sec:isotropic}, we introduce \emph{isotropic systems} and show that
they correspond to self-dual additive codes over $\GF(4)$.
It is known that isotropic systems can be represented by graphs. In Section~\ref{sec:graph} 
we define \emph{graph codes}. Theorem~\ref{th:graph} shows that every code can be represented
by a graph. This gives us a much smaller set of objects to work with. In Section~\ref{sec:lc}, 
we introduce \emph{local complementation}, and Theorem~\ref{th:lc} states that
two codes are equivalent if and only if the corresponding graphs are related via local complementations
and graph isomorphism. This implies that classifying codes up to equivalence is 
essentially the same as classifying orbits of graphs under local complementation.
We describe an algorithm for generating such graph orbits in Section~\ref{sec:classify}.
This algorithm was used to classify all codes of length up to 12. We show that Type~II codes
correspond to a special class of graphs and use this fact to classify all 
extremal Type~II codes of length 14. 
Finally, we determine that the smallest Type~I and Type~II 
codes with trivial automorphism group have length 9 and 12, respectively.
In Section~\ref{sec:concl}, we conclude and mention some other results.

\section{Isotropic Systems}\label{sec:isotropic}

We define a mapping $\phi : \GF(4) \to \GF(2)^2$ by $\phi(x) = (\tr(x\omega^2), \tr(x))$, i.e.,
$0 \mapsto (0,0)$, $1 \mapsto (1,0)$, $\omega \mapsto (0,1)$ and $\omega^2 \mapsto (1,1)$.
The reverse mapping $\phi^{-1} : \GF(2)^2 \to \GF(4)$ is given by $\phi^{-1}(a,b) = a + \omega b$.
Let $\boldsymbol{u} \in \GF(2)^{2n}$ be written as $\boldsymbol{u} = (\boldsymbol{a}|\boldsymbol{b}) =
(a_1,a_2,\ldots,a_n,b_1,b_2,\ldots,b_n)$. 
We extend the mapping $\phi : \GF(4)^n \to \GF(2)^{2n}$ by letting $\phi(\boldsymbol{v})
=(\boldsymbol{a}|\boldsymbol{b})$ where $\phi(v_i) = (a_i,b_i)$.
Likewise, we define $\phi^{-1} : \GF(2)^{2n} \to \GF(4)^n$ by $\phi^{-1}(\boldsymbol{a}|\boldsymbol{b})
= \boldsymbol{a} + \omega \boldsymbol{b}$.
We define the \emph{symplectic inner product}
of $(\boldsymbol{a}|\boldsymbol{b}), (\boldsymbol{a}'|\boldsymbol{b}') \in \GF(2)^{2n}$
as ${\langle(\boldsymbol{a}|\boldsymbol{b}), (\boldsymbol{a}'|\boldsymbol{b}')\rangle =
\boldsymbol{a} \cdot \boldsymbol{b}' + \boldsymbol{b} \cdot \boldsymbol{a}'}$.
A subset $\mathcal{I} \subset \GF(2)^{2n}$ is called \emph{totally isotropic} if 
$\langle\boldsymbol{u}, \boldsymbol{v}\rangle = 0$ for all 
$\boldsymbol{u}, \boldsymbol{v} \in \mathcal{I}$. 

\begin{defn}
A totally isotropic linear subspace of $\GF(2)^{2n}$ 
with dimension~$n$ defines an \emph{isotropic system}~\cite{bouchet1}.
An isotropic system can therefore be defined by the row space of a full rank $n \times 2n$ binary matrix $(A|B)$,
where $A B^{\text{T}} + B A^{\text{T}} = \boldsymbol{0}$.
\end{defn}

\begin{thm}
Every self-dual additive code over $\GF(4)$ can be uniquely represented as an isotropic system, and
every isotropic system can be uniquely represented as a self-dual additive code over $\GF(4)$.
\end{thm}
\begin{proof}
Let $\mathcal{C} \subset \GF(4)^n$ be a self-dual additive code.
Map $\mathcal{C}$ to $\mathcal{I} \subset \GF(2)^{2n}$ by mapping
each codeword $\boldsymbol{u} \in \mathcal{C}$ to
$\phi(\boldsymbol{u}) = (\boldsymbol{a}|\boldsymbol{b}) \in \GF(2)^{2n}$.
$\mathcal{I}$ must then be a linear subspace of $\GF(2)^{2n}$ with dimension~$n$.
$(\boldsymbol{a}|\boldsymbol{b}), (\boldsymbol{a}'|\boldsymbol{b}') \in \mathcal{I}$
are orthogonal with respect to the symplectic inner product if and only if 
$\phi^{-1}(\boldsymbol{a}|\boldsymbol{b}), \phi^{-1}(\boldsymbol{a}'|\boldsymbol{b}') \in \mathcal{C}$
are orthogonal with respect to the Hermitian trace inner product, because
\begin{eqnarray*}
&&\phi^{-1}(\boldsymbol{a}|\boldsymbol{b}) * \phi^{-1}(\boldsymbol{a}'|\boldsymbol{b}')\\
&=& \tr( \phi^{-1}(\boldsymbol{a}|\boldsymbol{b}) \cdot \overline{\phi^{-1}(\boldsymbol{a}'|\boldsymbol{b}')} )\\
&=& \tr( (\boldsymbol{a} + \omega \boldsymbol{b}) \cdot (\boldsymbol{a}' + \overline{\omega} \boldsymbol{b}'))\\
&=& (\boldsymbol{a} \cdot \boldsymbol{a}') \tr(1) + (\boldsymbol{a} \cdot \boldsymbol{b}') \tr(\overline{\omega}) 
+ (\boldsymbol{b} \cdot \boldsymbol{a}') \tr(\omega) + (\boldsymbol{b} \cdot \boldsymbol{b}') \tr(1) \\
&=& \boldsymbol{a} \cdot \boldsymbol{b}' + \boldsymbol{b} \cdot \boldsymbol{a}'.
\end{eqnarray*}
Since $\mathcal{C}$ is self-dual, $\boldsymbol{u} * \boldsymbol{v} = 0$ for all 
$\boldsymbol{u}, \boldsymbol{v} \in \mathcal{C}$, and $\mathcal{I}$ must therefore be totally isotropic.
It follows that $\mathcal{I}$ defines an isotropic system.
Likewise, the reverse mapping from an isotropic system to
a subset of $\GF(4)^n$ will always give a self-dual additive code over $\GF(4)$.
\end{proof}

\begin{exmp}\label{ex:code}
The row-space of $(A|B)$ defines an isotropic system, while $C = A + \omega B$ is a generator matrix of 
the $(6,2^6,4)$ Hexacode.
\[
(A|B) = 
\left(
\begin{array}{cccccc|cccccc}
1&0&0&1&0&0&0&0&0&0&1&1\\
0&0&0&0&1&1&1&0&0&1&1&1\\
0&1&0&0&1&0&0&0&0&1&0&1\\
0&0&0&1&0&1&0&1&0&1&1&1\\
0&0&1&0&0&1&0&0&0&1&1&0\\
0&0&0&1&1&0&0&0&1&1&1&1
\end{array}
\right)
\quad
C = 
\left(
\begin{array}{cccccc}
1&0&0&1&\omega&\omega\\
\omega&0&0&\omega&\omega^2&\omega^2\\
0&1&0&\omega&1&\omega\\
0&\omega&0&\omega^2&\omega&\omega^2\\
0&0&1&\omega&\omega&1\\
0&0&\omega&\omega^2&\omega^2&\omega
\end{array}
\right)
\]
\end{exmp}

\section{Graph Representation}\label{sec:graph}

A \emph{graph} is a pair $G=(V,E)$ where $V$ is a set of \emph{vertices},
and $E \subseteq V \times V$ is a set of \emph{edges}. A graph with $n$ vertices
can be represented by an $n \times n$ \emph{adjacency matrix} $\Gamma$, where
$\gamma_{ij} = 1$ if $\{i,j\} \in E$, and $\gamma_{ij} = 0$ otherwise.
We will only consider \emph{simple} \emph{undirected} graphs whose
adjacency matrices are symmetric with all diagonal elements being 0.
The \emph{neighbourhood} of $v \in V$, denoted $N_v \subset V$, is the set of 
vertices connected to $v$ by an edge. The number of vertices adjacent to $v$, $|N_v|$,
is called the \emph{degree} of $v$.
The \emph{induced subgraph} of $G$ on $W \subseteq V$ 
contains vertices $W$ and all edges from $E$ whose endpoints are both in $W$.
The \emph{complement} of $G$ is found by replacing $E$ with $V \times V - E$,
i.e., the edges in $E$ are changed to non-edges, and the non-edges to edges.
Two graphs $G=(V,E)$ and $G'=(V,E')$ are \emph{isomorphic} if and only if
there exists a permutation $\pi$ of $V$ such that $\{u,v\} \in E \iff \{\pi(u), \pi(v)\}
\in E'$.
A \emph{path} is a sequence of vertices, $(v_1,v_2,\ldots,v_i)$, such that
$\{v_1,v_2\}, \{v_2,v_3\},$ $\ldots, \{v_{i-1},v_{i}\} \in E$.
A graph is \emph{connected} if there is a path from any vertex to any other vertex in the graph.

\begin{defn}
A \emph{graph code} is an additive code over $\GF(4)$ that has a generator matrix of the form
$C = \Gamma + \omega I$, where $I$ is the identity matrix and $\Gamma$ is 
the adjacency matrix of a simple undirected graph.
\end{defn}

A graph code is always self-dual, since its generator matrix has full rank over $\GF(2)$
and $C \overline{C}^{\text{T}}$ only contains entries from $\GF(2)$ whose
traces must be zero.
This construction for self-dual additive codes over $\GF(4)$ has also
been used by Tonchev~\cite{tonchev}.

\begin{thm}\label{th:graph}
Every self-dual additive code over $\GF(4)$ is equivalent to a graph code.
\end{thm}
\begin{proof}
(This proof is due to Van~den~Nest, Dehaene, and De~Moor~\cite{nestthesis,nest}.)
We recall that the generator matrix of a self-dual additive code over $\GF(4)$
corresponds to an $n \times 2n$ binary matrix $(A|B)$, such that $C = A + \omega B$.
The row-space of $(A|B)$, denoted $\mathcal{I}$, defines an isotropic system.
We must prove that $\mathcal{I}$ is also generated by $(\Gamma|I)$, 
where $I$ is the identity matrix and $\Gamma$ is the adjacency matrix of a simple undirected graph.

The rows of $(A|B)$ can be replaced by any $n$ independent vectors from $\mathcal{I}$.
This basis change can be accomplished by 
$(A'|B') = M(A|B)$, where $M$ is an $n \times n$ invertible binary matrix.
If $B$ is invertible, the solution is simple, since $B^{-1} (A|B) = (\Gamma|I)$.
Note that $\Gamma$ will always be a symmetric matrix, since $\Gamma I^{\text{T}} + I \Gamma^{\text{T}} = 0$.
If the $i$th diagonal element of $\Gamma$ is 1, it can be set to 0 by
conjugating column~$i$ of $\Gamma + \omega I$.

In the case where $B$ has rank $k<n$, we can perform a basis change to get
\[
(A'|B') = 
\left(\begin{array}{c|c}A_1 & B_1 \\ A_2 & \boldsymbol{0} \end{array}\right),
\]
where $B_1$ is a $k \times n$ matrix with full rank, and 
$A_1$ also has size $k \times n$.
Since the row-space of $(A'|B')$ is totally isotropic, and $B'$ contains an all-zero
row, it must be true that $A_2 B_1^{\text{T}} = \boldsymbol{0}$.
$A_2$ must have full rank, and the row space
of $B_1$ must be the orthogonal complement of the row space of $A_2$.

We assume that $B_1 = (B_{11} | B_{12})$ where $B_{11}$ is a $k \times k$ invertible
matrix. We also write $A_2 = (A_{21} | A_{22})$ where $A_{22}$ has size
$(n-k) \times (n-k)$. 
Assume that there exists an $\boldsymbol{x} \in \GF(2)^{n-k}$ such that $A_{22} \boldsymbol{x}^{\text{T}} = 0$.
Then the vector $\boldsymbol{v} = (0, \ldots, 0, \boldsymbol{x})$ of length~$n$
satisfies $A_2 \boldsymbol{v}^{\text{T}} = 0$.
Since the row space of $B_1$ is the orthogonal complement of the row space of $A_2$,
we can write $\boldsymbol{v} = \boldsymbol{y} B_1$ for some $\boldsymbol{y} \in \GF(2)^k$.
We see that $\boldsymbol{y} B_{11}  = 0$, and since $B_{11}$ has full rank, it must therefore be true that
$\boldsymbol{y} = 0$. This means that $\boldsymbol{x} = 0$, which proves that
$A_{22}$ is an invertible matrix.

Interchanging column~$i$ of $A'$ and column~$i$ of $B'$ corresponds to multiplication by 
$\omega^2$ followed by conjugation of the $i$th column of $A' + \omega B'$.
We can therefore swap the $i$th columns of $A'$ and $B'$ for $k < i \le n$ to get
$(A''|B'')$. Since $B_{11}$ and $A_{22}$ are invertible, $B''$ must also be an invertible matrix.
We then find $B''^{-1} (A''|B'') = (\Gamma|I)$, and set all diagonal elements of $\Gamma$ to 0.
\end{proof}

\begin{figure}
 \centering
 {\hfill
 \subfloat[]{\includegraphics[width=.40\linewidth]{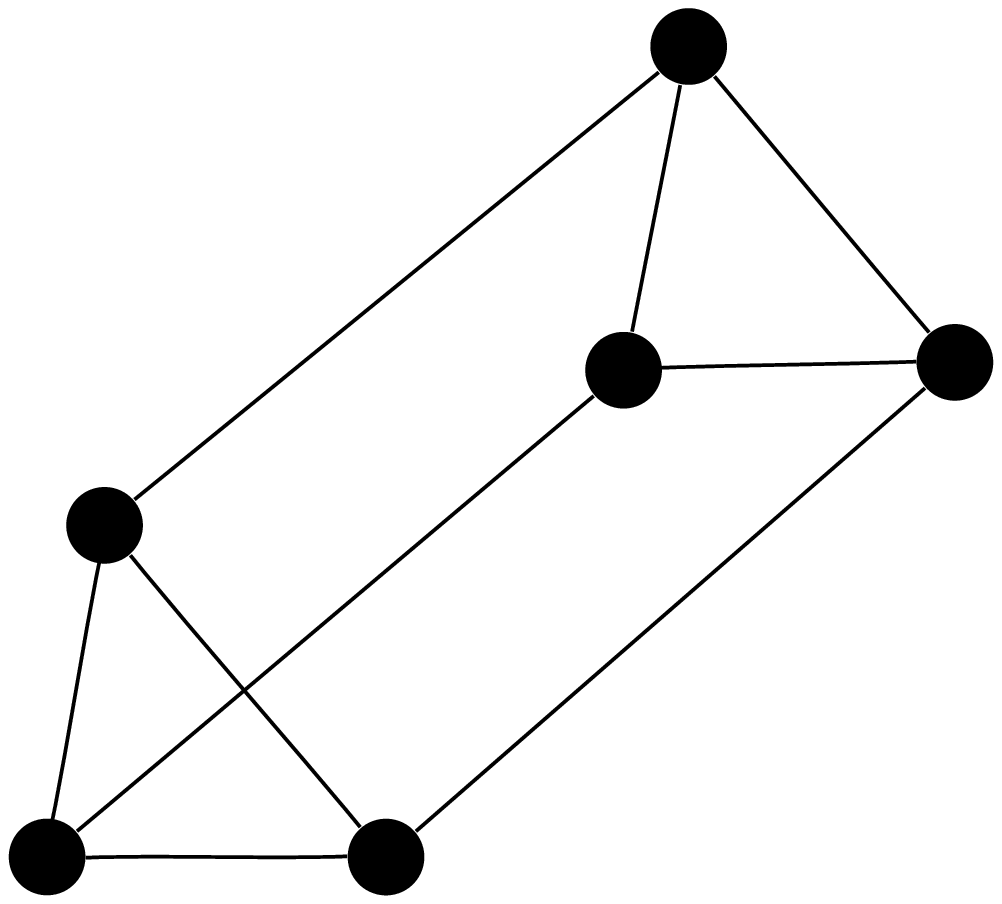}\label{fig:doubleclique}}
 \hfill
 \subfloat[]{\includegraphics[width=.40\linewidth]{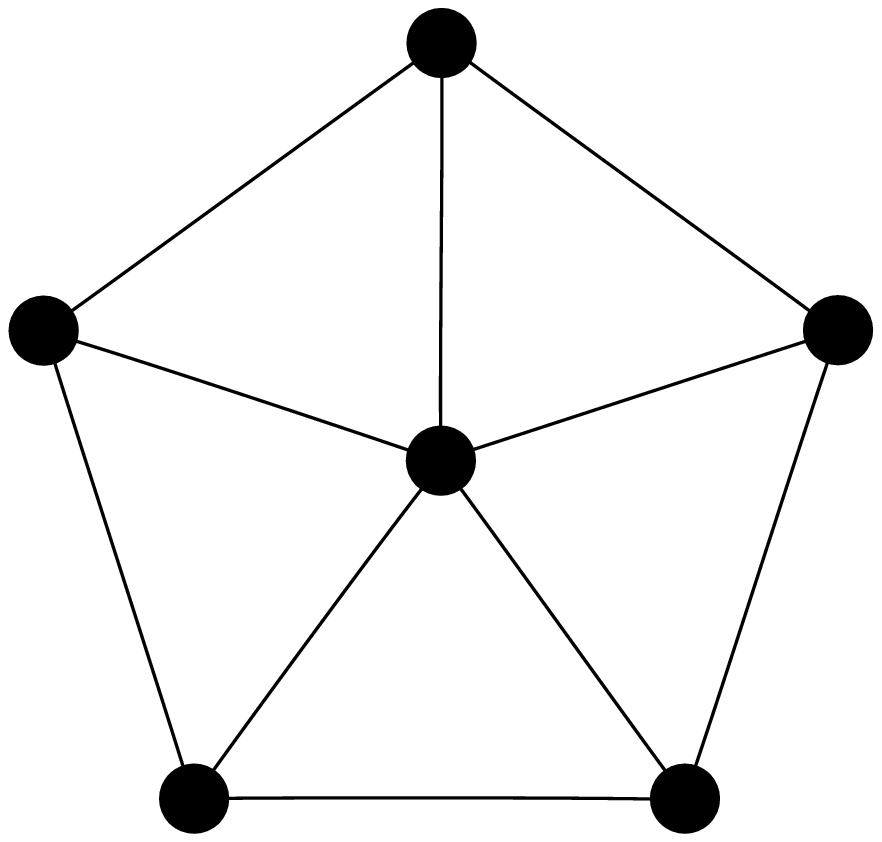}\label{fig:wheel}}
 \hfill}
 \caption{Two Graph Representations of the Hexacode}
\end{figure}

\begin{exmp}\label{ex:graphcode}
Let $C=A + \omega B$ be the generator matrix of the $(6,2^6,4)$ Hexacode
given in Example~\ref{ex:code}.
By the method described in the proof of Theorem~\ref{th:graph},
we find $C' = \Gamma + \omega I$, which generates an equivalent graph code.
$\Gamma$ is the adjacency matrix of the graph shown in Fig.~\ref{fig:wheel}.
\[
(A|B) = 
\left(
\begin{array}{cccccc|cccccc}
1&0&0&1&0&0&0&0&0&0&1&1\\
0&0&0&0&1&1&1&0&0&1&1&1\\
0&1&0&0&1&0&0&0&0&1&0&1\\
0&0&0&1&0&1&0&1&0&1&1&1\\
0&0&1&0&0&1&0&0&0&1&1&0\\
0&0&0&1&1&0&0&0&1&1&1&1
\end{array}
\right)
\quad
(\Gamma|I) = 
\left(
\begin{array}{cccccc|cccccc}
0&0&1&0&1&1 &1&0&0&0&0&0\\
0&0&1&1&0&1 &0&1&0&0&0&0\\
1&1&0&0&0&1 &0&0&1&0&0&0\\
0&1&0&0&1&1 &0&0&0&1&0&0\\
1&0&0&1&0&1 &0&0&0&0&1&0\\
1&1&1&1&1&0 &0&0&0&0&0&1
\end{array}
\right)
\]
\[
C' = 
\left(
\begin{array}{cccccc}
\omega&0&1&0&1&1\\
0&\omega&1&1&0&1\\
1&1&\omega&0&0&1\\
0&1&0&\omega&1&1\\
1&0&0&1&\omega&1\\
1&1&1&1&1&\omega
\end{array}
\right)
\]
\end{exmp}

Theorem~\ref{th:graph} was first proved by Bouchet~\cite{bouchet4} in the 
context of isotropic systems, and later by Schlingemann~\cite{schlingemann2}
in terms of \emph{quantum stabilizer states}.
Proofs of Theorem~\ref{th:graph} have also been given by 
Schlingemann and Werner~\cite{schlingemann}, by
Grassl, Klappenecker, and Rötteler~\cite{grasslpaper}, by Glynn et~al.~\cite{glynnbook,glynngraph},
and by Van~den~Nest et~al.~\cite{nestthesis,nest}.

Swapping vertex~$i$ and vertex $j$ of a graph with adjacency matrix $\Gamma$ can be accomplished by
exchanging column~$i$ and column~$j$ of $\Gamma$ and then exchanging row~$i$ and row~$j$ of $\Gamma$.
We call the resulting matrix $\Gamma'$. Exactly the same column and row operations map
$\Gamma + \omega I$ to $\Gamma' + \omega I$. These matrices generate equivalent codes.
It follows that two codes are equivalent if their corresponding graphs are isomorphic.

We have seen that every graph represents a self-dual additive code over $\GF(4)$,
and that every self-dual additive code over $\GF(4)$ can be represented by a graph.
It follows that we can, without loss of generality,
restrict our study to codes with generator matrices of the form $\Gamma + \omega I$,
where $\Gamma$ are adjacency matrices of unlabeled simple undirected graphs.

\section{Local Complementation}\label{sec:lc}

\begin{defn}
Given a graph $G=(V,E)$ and a vertex $v \in V$,
let $N_v \subset V$ be the neighbourhood of $v$.
\emph{Local complementation} (LC) on $v$ transforms $G$ into $G^v$. 
To obtain $G^v$, we replace the induced subgraph of $G$ on $N_v$ by its complement.
It is easy to verify that $(G^v)^v = G$.
\end{defn}

\begin{figure}
 \centering
 \subfloat[The Graph $G$]
 {\hspace{5pt}\includegraphics[width=.25\linewidth]{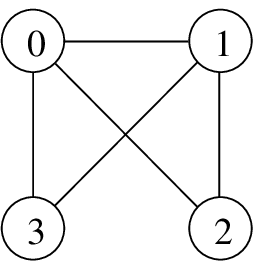}\hspace{5pt}\label{fig:lcexample1}}
 \quad
 \subfloat[The LC Image $G^0$]
 {\hspace{5pt}\includegraphics[width=.25\linewidth]{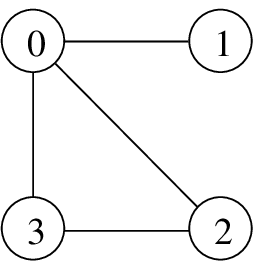}\hspace{5pt}\label{fig:lcexample2}}
 \caption{Example of Local Complementation}\label{fig:lcexample}
\end{figure}

\begin{exmp}
We will perform local complementation on vertex 0 of the
graph $G$, shown in Fig.~\ref{fig:lcexample1}.
We see that the neighbourhood of 0 is $N_0 = \{1,2,3\}$
and that the induced subgraph on the neighbourhood has edges $\{1,2\}$ and $\{1,3\}$. 
The complement of this subgraph contains the single edge $\{2,3\}$.
The resulting LC image, $G^0$, is seen in Fig.~\ref{fig:lcexample2}.
\end{exmp}

\begin{exmp}
Consider the
graph shown in Fig.~\ref{fig:doubleclique}, whose corresponding
graph code is the Hexacode. An LC operation on any vertex of this
graph produces the graph shown in Fig.~\ref{fig:wheel}.
An LC operation on the vertex in the centre of the graph shown in
Fig.~\ref{fig:wheel} gives the same graph, up to isomorphism.
LC operations on any of the other five vertices produces
the graph shown in Fig.~\ref{fig:doubleclique}.
\end{exmp}

\begin{thm}\label{thm:lcmatrix}
Let $\Gamma$ be the adjacency matrix of the graph $G=(V,E)$, and
$\Gamma^v$ be the adjacency matrix of $G^v$, for any $v \in V$.
The codes generated by $C = \Gamma + \omega I$ and
$C' = \Gamma^v + \omega I$ are equivalent.
\end{thm}
\begin{proof}
We show that $C$ can be transformed into $C'$ by using only operations that
map a code to an equivalent code.
Each row and each column of $C$ correspond to a vertex in $V$.
Let $N_v$ denote the neighbourhood of $v$.
For all  $i \in N_v$, add row~$v$ of $C$ to row~$i$ of $C$.
Multiply column~$v$ of $C$ by $\omega$ and then conjugate the
same column. Finally, conjugate column~$i$ of $C$, for all $i \in N_v$.
The resulting matrix is $C'$.
\end{proof}

\begin{thm}\label{th:lc}
Two self-dual additive codes over $\GF(4)$, $\mathcal{C}$ and $\mathcal{C}'$,
with graph representations $G$ and $G'$, are equivalent if and only if
there is a finite sequence of not necessarily distinct vertices
$(v_1, v_2, \ldots, v_i)$, such that $(((G^{v_1})^{v_2})^{\cdots})^{v_{i}}$ is
isomorphic to $G'$.
\end{thm}
\begin{proof}[Sketch of proof]
Let $\Gamma$ be the adjacency matrix of $G$, and let $\mathcal{C}_G$ be the code generated by 
$\Gamma + \omega I$. Likewise, let $\Gamma'$ be the adjacency matrix of $G'$, and let $\mathcal{C}_G'$
be the code generated by $\Gamma' + \omega I$.
If the codewords of $\mathcal{C}$ are mapped onto the codewords of $\mathcal{C}'$
by one of the $6^n n!$ combinations of coordinate permutations, coordinate scalings, 
and coordinate conjugations, then there must also be a transformation from this set that maps the codewords of 
$\mathcal{C}_G$ onto the codewords of $\mathcal{C}'_G$.
Consequently, we only need to consider those transformations that map a graph code
to another graph code. The codes obtained by the $n!$ possible permutations of coordinates correspond
to graph isomorphisms.

Let $C = \Gamma + \omega I$ be transformed into $C' = A + \omega B$ by coordinate scalings and conjugations.
Then $C'$ is a graph
code if and only if $B$ is invertible and all diagonal elements of $B^{-1} A$ are zero.
It is easy to verify that conjugation of column~$i$ of $C'$ has no effect on $B$, but flips
the value of the $i$th diagonal element of $B^{-1} A$.
Given a combination of column scalings on $C$ such that the resulting $B$ is invertible,
there must therefore be a unique combination of column conjugations on $C$ such that the resulting
$B^{-1} A$ has zero diagonal.
We must therefore show that any combination of column scalings on $C$ that give an invertible $B$
can be performed by a sequence of LC operations on $G$.

Multiplying column~$i$ of $C$ by $\omega^2$ replaces column~$i$ of $I$ with column~$i$ of $\Gamma$.
Multiplying column~$i$ of $C$ by $\omega$ adds column~$i$ of $\Gamma$ to column~$i$ of $I$.
It is then possible to show which of the $3^n$ possible scalings do not give an invertible $B$.
A vertex $v$ of $G$ corresponds to a column of $\Gamma$. The neighbourhood of $v$,
$N_v$, corresponds to a set of columns of $\Gamma$.
We know from Theorem~\ref{thm:lcmatrix} that an LC operation on vertex~$i$ of $G$ corresponds
to a scaling of column~$i$ of $C$ by $\omega$ followed by conjugation of column~$i$ and all columns in $N_i$.
Observe that conjugating a coordinate followed by a scaling by $\omega$ is equivalent to
scaling by $\omega^2$ followed by conjugation.
Note in particular that the local complementations $((G^i)^j)^i$, where $i$ and $j$ are adjacent vertices,
are equivalent to scaling both column~$i$ and column~$j$ of $C$ by $\omega^2$.
It can be verified that any combination of column scalings that
map a graph code to a graph code can be implemented as a sequence of LC operations.
The exact algorithm for finding this sequence of LC operations is quite involved,
and we refer to the proof by Van~den~Nest et~al.~\cite{nestthesis,nest} for the details.
\end{proof}

Bouchet~\cite{bouchet4} first proved Theorem~\ref{th:lc} in terms of isotropic systems.
The same result was discovered 
by Van~den~Nest et~al.~\cite{nestthesis,nest} in terms of quantum stabilizer states, and by
Glynn~et~al.~\cite{glynnbook,glynngraph} using finite geometry.

\section{Classification}\label{sec:classify}

\begin{defn}
The \emph{LC orbit} of a graph $G$ is the set of all unlabeled graphs that
can be obtained by performing any sequence of LC operations on~$G$.
\end{defn}

It follows from Theorem~\ref{th:lc} that two self-dual additive codes over $\GF(4)$
are equivalent if and only if their graph representations are in the same LC orbit.
As an example, the two graphs shown in Fig.~\ref{fig:doubleclique} and Fig.~\ref{fig:wheel}
make up a complete LC orbit, and are thus the only possible graph
representations of the Hexacode.
The LC orbit of a graph can easily be generated by a recursive algorithm.
We have used the program \emph{nauty}~\cite{nauty} to check for graph isomorphism.

Let $\boldsymbol{G}_n$ be the set of all unlabeled simple undirected connected 
graphs on $n$ vertices.
Connected graphs correspond to \emph{indecomposable} codes.
A code is decomposable if it can be written as the \emph{direct sum} of two smaller codes.
For example, let $\mathcal{C}$ be an $(n,2^n,d)$ code and $\mathcal{C}'$ an $(n',2^{n'},d')$ code. The
direct sum, $\mathcal{C} \oplus \mathcal{C}' = \{u||v \mid u \in \mathcal{C}, v \in \mathcal{C}'\}$,
where $||$ means concatenation, is an $({n+n'},2^{n+n'},\min\{d,d'\})$ code.
It follows that all decomposable codes of length~$n$ can be classified easily once
all indecomposable codes of length less than $n$ are known.

The set of all distinct LC orbits of connected graphs on $n$ vertices is a partitioning of
$\boldsymbol{G}_n$ into $i_n$ disjoint sets.
$i_n$ is also the number of indecomposable self-dual additive codes over $\GF(4)$ of length~$n$,
up to equivalence. Let $\boldsymbol{L}_n$ be a set containing one representative from each
LC orbit of connected graphs on $n$ vertices. We have devised several algorithms~\cite{mscthesis} 
for finding such sets of representatives.
The simplest approach is to start with the set $\boldsymbol{G}_n$ and generate LC orbits
of its members until we have a partitioning of $\boldsymbol{G}_n$.
The following more efficient technique was described by Glynn et~al.~\cite{glynnbook}.
Let the $2^n-1$ \emph{extensions} of a graph on $n$ vertices be formed by
adding a new vertex and joining it to all possible combinations of at least one of the old vertices.
The set $\boldsymbol{E}_{n}$, containing $i_{n-1} (2^{n-1}-1)$ graphs,
is formed by making all possible extensions of all graphs in $\boldsymbol{L}_{n-1}$.

\begin{thm}\label{thm:extend}
$\boldsymbol{L}_n \subset \boldsymbol{E}_{n}$, i.e.,
the set $\boldsymbol{E}_{n}$ will contain at least one representative from each
LC orbit of connected graphs on $n$ vertices.
\end{thm}
\begin{proof}
Let $G=(V,E) \in \boldsymbol{G}_n$, and choose any subset $W \subset V$ of $n-1$ vertices.
By doing LC operations on vertices in $W$, we can transform the induced subgraph of $G$ on $W$ into 
one of the graphs in $\boldsymbol{L}_{n-1}$ that were extended when $\boldsymbol{E}_{n}$ was constructed.
It follows that for all $G \in \boldsymbol{G}_{n}$, some graph in the LC orbit of $G$ must be 
part of $\boldsymbol{E}_{n}$.
\end{proof}

The set $\boldsymbol{E}_{n}$ will be much smaller than $\boldsymbol{G}_n$, so it will be more
efficient to search for a set of LC orbit representatives within $\boldsymbol{E}_{n}$.
It is also desirable to partition the set $\boldsymbol{E}_{n}$ such that graphs from two different
partitions are guaranteed to belong to different LC orbits. We can then consider each partition independently,
which reduces the amount of memory required and allow for parallel processing.
To do this, we must have some property that is invariant over the LC orbit and
that can be calculated quickly.

The special form of the generator matrix of a graph code makes it easier to find 
the number of codewords of weight $i<n$. If $\mathcal{C}$ is generated by $C = \Gamma + \omega I$,
then any codeword formed by adding $i$ rows of $C$ must have weight at least $i$.
This means that we can find the \emph{partial weight distribution} of $\mathcal{C}$,
$(A_0, A_1, \ldots, A_j)$, for some $j<n$, by only considering codewords
formed by adding $j$ or fewer rows of $C$.
We calculate the partial weight distribution, for a suitable choice of $j$, of
all codes corresponding to graphs in $\boldsymbol{E}_{n}$.
Codes with different partial weight distribution can never be equivalent,
so we partition $\boldsymbol{E}_{n}$ such that graphs corresponding to codes
with the same partial weight distribution are always in the same partition.

Using the described techniques, and a parallel cluster computer,
we were able to classify all self-dual additive codes
over $\GF(4)$ of length up to 12. The results have been verified by checking that 
the sizes of all LC orbits add up to the number of graphs in $\boldsymbol{G}_n$.
The sizes of the automorphism groups of all codes have also been calculated, and
it has been verified that that the mass formulas defined by
Eq.~\eqref{eq:mass} and Eq.~\eqref{eq:mass2} give the correct values.
Table~\ref{tab:orbitsindecomp} gives the values of $i_n$, the number of distinct LC orbits
of connected graphs on $n$ vertices, which is also the number of inequivalent indecomposable 
codes of length~$n$. The table also gives the values of $i_n^{\text{II}}$, the number of 
indecomposable Type~II codes. The total number of inequivalent codes of length~$n$, $t_n$,
and the total number of Type~II codes of length~$n$, $t_n^{\text{II}}$, are shown in 
Table~\ref{tab:orbitsall}. The numbers $t_n$ are easily derived from the numbers $i_n$ by using
the \emph{Euler transform}~\cite{sloane2},
\begin{eqnarray*}
c_n &=& \sum_{d|n} d i_d\\
t_1 &=& c_1\\
t_n &=& \frac{1}{n}\left( c_n + \sum_{k=1}^{n-1} c_k t_{n-k} \right).
\end{eqnarray*}
The numbers $t_n^{\text{II}}$ are similarly derived from $i_n^{\text{II}}$.
The values of $i_n$ and $t_n$ can be found as sequences A090899
and A094927 in \emph{The On-Line Encyclopedia of Integer Sequences}~\cite{sequences}.
Table~\ref{tab:distances_ind} and Table~\ref{tab:distances} list by minimum distance the 
numbers of indecomposable codes and the total numbers of codes.\footnote{Note that some 
authors~\cite{gaborit,huffman} give $3$ as the total number of self-dual $(7,2^7,3)$ codes.
The correct number is 4.} Table~\ref{tab:type2_ind} 
and Table~\ref{tab:type2} similarly list the numbers of Type~II codes by minimum distance.
The numbers of Type~I codes can be obtained by subtracting the numbers of Type~II codes
from the total numbers.
The number of distinct weight enumerators of all codes of length~$n$ and minimum distance~$d$
can be found in Table~\ref{tab:wd}. There are obviously too many codes to give a 
complete list here, but a database containing one representative from each equivalence class, 
with information about weight enumerators, automorphism groups, etc., is available on-line~\cite{database}.

\begin{table}
\centering
\caption{Number of Indecomposable ($i_n$) and Indecomposable Type~II ($i_n^{\text{II}}$) Codes of Length~$n$}
\label{tab:orbitsindecomp}
\begin{tabular}{ccccccccccccc}
\toprule
$n$ & 1 & 2 & 3 & 4 & 5 & 6 & 7 & 8 & 9 & 10 & 11 & 12 \\
\midrule
$i_n$ & 1 & 1 & 1 & 2 &  4 & 11 & 26 & 101 & 440 & 3,132 & 40,457 & 1,274,068 \\
$i_n^{\text{II}}$ &   & 1 &   & 1 &   & 4 &   & 14 &   & 103 &   & 2,926\\
\bottomrule
\end{tabular}
\end{table}

\begin{table}
\centering
\caption{Total Number ($t_n$) and Number of Type~II ($t_n^{\text{II}}$) Codes of Length~$n$}
\label{tab:orbitsall}
\begin{tabular}{ccccccccccccc}
\toprule
$n$ & 1 & 2 & 3 & 4 & 5 & 6 & 7 & 8 & 9 & 10 & 11 & 12 \\
\midrule
$t_n$ & 1 & 2 & 3 & 6 & 11 & 26 & 59 & 182 & 675 & 3,990 & 45,144 & 1,323,363 \\
$t_n^{\text{II}}$ &   & 1 &   & 2 &   & 6 &   & 21 &   & 128 &   & 3,079 \\
\bottomrule
\end{tabular}
\end{table}

\begin{table}
\centering
\caption{Number of Indecomposable Codes of Length~$n$ and Minimum Distance~$d$}
\label{tab:distances_ind}
\begin{tabular}{crrrrrrrrrrr}
\toprule
$d \backslash n$ & 
\multicolumn{1}{c}{2} & \multicolumn{1}{c}{3} & \multicolumn{1}{c}{4} & \multicolumn{1}{c}{5} &
\multicolumn{1}{c}{6} & \multicolumn{1}{c}{7} & \multicolumn{1}{c}{8} & \multicolumn{1}{c}{9} &
\multicolumn{1}{c}{10} & \multicolumn{1}{c}{11} & \multicolumn{1}{c}{12} \\
\midrule
2     & 1 & 1 & 2 &  3 &  9 & 22 &  85 & 363 & 2,436 & 26,750 &   611,036 \\
3     &   &   &   &  1 &  1 &  4 &  11 &  69 &   576 & 11,200 &   467,513 \\
4     &   &   &   &    &  1 &    &   5 &   8 &   120 &  2,506 &   195,455 \\
5     &   &   &   &    &    &    &     &     &       &      1 &        63 \\
6     &   &   &   &    &    &    &     &     &       &        &         1 \\
\midrule
All   & 1 & 1 & 2 &  4 & 11 & 26 & 101 & 440 & 3,132 & 40,457 & 1,274,068 \\
\bottomrule
\end{tabular}
\end{table}

\begin{table}
\centering
\caption{Total Number of Codes of Length~$n$ and Minimum Distance~$d$}
\label{tab:distances}
\begin{tabular}{crrrrrrrrrrrr}
\toprule
$d \backslash n$ & 
\multicolumn{1}{c}{1} &
\multicolumn{1}{c}{2} & \multicolumn{1}{c}{3} & \multicolumn{1}{c}{4} & \multicolumn{1}{c}{5} &
\multicolumn{1}{c}{6} & \multicolumn{1}{c}{7} & \multicolumn{1}{c}{8} & \multicolumn{1}{c}{9} &
\multicolumn{1}{c}{10} & \multicolumn{1}{c}{11} & \multicolumn{1}{c}{12} \\
\midrule
1     & 1 & 1 & 2 & 3 &  6 & 11 & 26 &  59 & 182 &   675 &  3,990 &    45,144 \\
2     &   & 1 & 1 & 3 &  4 & 13 & 29 & 107 & 416 & 2,618 & 27,445 &   615,180 \\
3     &   &   &   &   &  1 &  1 &  4 &  11 &  69 &   577 & 11,202 &   467,519 \\
4     &   &   &   &   &    &  1 &    &   5 &   8 &   120 &  2,506 &   195,456 \\
5     &   &   &   &   &    &    &    &     &     &       &      1 &        63 \\
6     &   &   &   &   &    &    &    &     &     &       &        &         1 \\
\midrule
All   & 1 & 2 & 3 & 6 & 11 & 26 & 59 & 182 & 675 & 3,990 & 45,144 & 1,323,363 \\
\bottomrule
\end{tabular}
\end{table}

\begin{table}
\centering
\caption{Number of Indecomposable Type~II Codes of Length~$n$ and Minimum Distance~$d$}
\label{tab:type2_ind}
\begin{tabular}{crrrrrrr}
\toprule
$d \backslash n$ 
& \multicolumn{1}{c}{2} & \multicolumn{1}{c}{4} & \multicolumn{1}{c}{6} 
 & \multicolumn{1}{c}{8} & \multicolumn{1}{c}{10} & \multicolumn{1}{c}{12} & \multicolumn{1}{c}{14}\\
\midrule
2     & 1 & 1 & 3 & 11 &  84 & 2,133 & \multicolumn{1}{c}{?} \\
4     &   &   & 1 &  3 &  19 &   792 & \multicolumn{1}{c}{?} \\
6     &   &   &   &    &     &     1 & 1,020 \\
\midrule
Total & 1 & 1 & 4 & 14 & 103 & 2,926 & \multicolumn{1}{c}{?} \\
\bottomrule
\end{tabular}
\end{table}

\begin{table}
\centering
\caption{Total Number of Type~II Codes of Length~$n$ and Minimum Distance~$d$}
\label{tab:type2}
\begin{tabular}{crrrrrrr}
\toprule
$d \backslash n$ 
& \multicolumn{1}{c}{2} & \multicolumn{1}{c}{4} & \multicolumn{1}{c}{6} 
 & \multicolumn{1}{c}{8} & \multicolumn{1}{c}{10} & \multicolumn{1}{c}{12} & \multicolumn{1}{c}{14}\\
\midrule
2     & 1 & 2 & 5 & 18 & 109 & 2,285 & \multicolumn{1}{c}{?}\\
4     &   &   & 1 &  3 &  19 &   793 & \multicolumn{1}{c}{?}\\
6     &   &   &   &    &     &     1 & 1,020\\
\midrule
Total & 1 & 2 & 6 & 21 & 128 & 3,079 & $\ge$ 1,727,942\\
\bottomrule
\end{tabular}
\end{table}

\begin{table}
\centering
\caption{Number of Distinct Weight Enumerators of All Codes of Length~$n$ and Minimum Distance~$d$}
\label{tab:wd}
\begin{tabular}{crrrrrrrrrrrr}
\toprule
$d \backslash n$ & 
\multicolumn{1}{c}{1} &
\multicolumn{1}{c}{2} & \multicolumn{1}{c}{3} & \multicolumn{1}{c}{4} & \multicolumn{1}{c}{5} &
\multicolumn{1}{c}{6} & \multicolumn{1}{c}{7} & \multicolumn{1}{c}{8} & \multicolumn{1}{c}{9} &
\multicolumn{1}{c}{10} & \multicolumn{1}{c}{11} & \multicolumn{1}{c}{12} \\
\midrule
1     & 1 & 1 & 2 & 3 &  5 & 10 & 23 &  46 & 116 &   320 &    909 &     3,312 \\
2     &   & 1 & 1 & 2 &  4 & 11 & 21 &  64 & 187 &   549 &  2,249 &    11,419 \\
3     &   &   &   &   &  1 &  1 &  2 &   4 &  15 &    33 &    125 &       625 \\
4     &   &   &   &   &    &  1 &    &   2 &   2 &     7 &     28 &       178 \\
5     &   &   &   &   &    &    &    &     &     &       &      1 &         2 \\
6     &   &   &   &   &    &    &    &     &     &       &        &         1 \\
\midrule
All   & 1 & 2 & 3 & 5 & 10 & 23 & 46 & 116 & 320 &   909 &  3,312 &    15,537 \\
\bottomrule
\end{tabular}
\end{table}

Our results give a complete classification of the extremal Type~I $(10,2^{10},4)$ and $(12,2^{12},5)$ codes. 
These classifications were previously unknown.
The 101 extremal Type~I $(10,2^{10},4)$ codes have 6 distinct weight enumerators,
{\footnotesize\begin{eqnarray*}
W_{10,1}(x,y) &=& x^{10} + 10 x^6 y^4 + 72 x^5 y^5 + 160 x^4 y^6 + 240 x^3 y^7 + 285 x^2 y^8 + 200 x y^9 + 56 y^{10},\\
W_{10,2}(x,y) &=& x^{10} + 14 x^6 y^4 + 64 x^5 y^5 + 156 x^4 y^6 + 256 x^3 y^7 + 281 x^2 y^8 + 192 x y^9 + 60 y^{10},\\
W_{10,3}(x,y) &=& x^{10} + 18 x^6 y^4 + 56 x^5 y^5 + 152 x^4 y^6 + 272 x^3 y^7 + 277 x^2 y^8 + 184 x y^9 + 64 y^{10},\\
W_{10,4}(x,y) &=& x^{10} + 22 x^6 y^4 + 48 x^5 y^5 + 148 x^4 y^6 + 288 x^3 y^7 + 273 x^2 y^8 + 176 x y^9 + 68 y^{10},\\
W_{10,5}(x,y) &=& x^{10} + 26 x^6 y^4 + 40 x^5 y^5 + 144 x^4 y^6 + 304 x^3 y^7 + 269 x^2 y^8 + 168 x y^9 + 72 y^{10},\\
W_{10,6}(x,y) &=& x^{10} + 30 x^6 y^4 + 32 x^5 y^5 + 140 x^4 y^6 + 320 x^3 y^7 + 265 x^2 y^8 + 160 x y^9 + 76 y^{10}.
\end{eqnarray*}}
Table~\ref{tab:10codes} lists the number of such codes by weight enumerator and automorphism group size.
The 63 extremal Type~I $(12,2^{12},5)$ codes have 2 distinct weight enumerators,
\begin{gather*}
\begin{split}
W_{12,1}(x,y) &= x^{12} + 40x^7 y^5 + 212 x^6 y^6 + 424 x^5 y^7 + 725 x^4 y^8 + 1080 x^3 y^9 + \\
             &\quad 980 x^2 y^{10} + 504 x y^{11} + 130 y^{12},
\end{split}\\
\begin{split}
W_{12,2}(x,y) &= x^{12} + 48x^7 y^5 + 188 x^6 y^6 + 432 x^5 y^7 + 765 x^4 y^8 + 1040 x^3 y^9 + \\
             &\quad 972 x^2 y^{10} + 528 x y^{11} + 122 y^{12}.
\end{split}
\end{gather*}
Table~\ref{tab:12codes} lists the number of such codes by weight enumerator and automorphism group size.

\begin{table}
\centering
\caption{Number of Extremal Type~I $(10,2^{10},4)$ Codes with Weight Enumerator~$w$ and Automorphism Group of Size~$a$}
\label{tab:10codes}
\begin{tabular}{cccccccc}
\toprule
$a \backslash w$ & $W_{10,1}$ & $W_{10,2}$ & $W_{10,3}$ & $W_{10,4}$ & $W_{10,5}$ & $W_{10,6}$ & All\\
\midrule
1    &   & 3 &   &   &   &   & 3\\
2    & 2 & 9 & 7 & 2 &   &   & 20\\
4    & 5 & 9 & 7 & 1 &   &   & 22\\
6    & 1 &   &   & 1 &   &   & 2\\
8    & 1 & 4 & 3 &   & 1 &   & 9\\
12   &   & 1 &   &   &   &   & 1\\
16   & 1 & 1 & 6 & 5 & 3 &   & 16\\
32   & 2 & 2 & 2 & 1 & 2 &   & 9\\
40   & 1 &   &   &   &   &   & 1\\
48   & 1 &   &   & 3 &   &   & 4\\
64   &   &   & 2 &   &   &   & 2\\
128  &   & 2 &   &   &   &   & 2\\
192  &   & 1 & 2 &   & 1 &   & 4\\
256  &   &   &   &   &   & 2 & 2\\
320  & 1 &   &   &   &   & 1 & 2\\
384  &   &   &   &   &   & 1 & 1\\
3840 &   &   &   &   &   & 1 & 1\\
\midrule
All  &15 &32 &29 &13 & 7 & 5 & 101\\
\bottomrule
\end{tabular}
\end{table}

\begin{table}
\centering
\caption{Number of Extremal Type~I $(12,2^{12},5)$ Codes with Weight Enumerator~$w$ and Automorphism Group of Size~$a$}
\label{tab:12codes}
\begin{tabular}{cccc}
\toprule
$a \backslash w$ & $W_{12,1}$ & $W_{12,2}$ & All\\
\midrule
1    &   & 25& 25\\
2    &   & 23& 23\\
3    &   &  1& 1\\
4    & 3 &  4& 7\\
6    & 1 &  3& 4\\
8    &   &  2& 2\\
24   &   &  1& 1\\
\midrule
All  & 4 & 59&63\\
\bottomrule
\end{tabular}
\end{table}

By observing that graphs corresponding to Type~II codes have a special property, we are able
to extend our classification to all the 1,020 extremal Type~II $(14,2^{14},6)$ codes. It was previously shown by
Bachoc and Gaborit~\cite{bachoc} that there are at least 490 such codes.

\begin{thm}
Let $\Gamma$ be the adjacency matrix of the graph $G$.
The code $\mathcal{C}$ generated by $C = \Gamma + \omega I$ is of Type~II if and only if
$G$ is \emph{anti-Eulerian}, i.e., if all its vertices have odd degree.
\end{thm}
\begin{proof}
If $\mathcal{C}$ is of Type~II, then every row of $C$ must have even weight. It follows that
every row of $\Gamma$ must have odd weight, and therefore correspond to an anti-Eulerian graph.
Conversely, if all rows of $C$ have even weight, $\mathcal{C}$ must be of Type~II, since
the codeword formed by adding any subset of these rows must also have even weight.
This follows from the fact that for any two codewords of a self-dual code, there must be an
even number of coordinates where the codewords have different non-zero values.
\end{proof}

An anti-Eulerian graph is the complement of an \emph{Eulerian graph}, i.e., a graph where all vertices have even degree.
It is easy to show that all anti-Eulerian graphs must have an even number of vertices,
and it follows that all Type~II codes must have even length.
To classify Type~II codes of length 14, we proceed as follows.
We take the set $\boldsymbol{L}_{12}$ containing 1,274,068 LC orbit representatives of graphs on 12 vertices.
All these graphs are then extended, but in a slightly different way than earlier.
To each graph we add one vertex and join it to all possible combinations of at least one of the old vertices.
To each obtained graph we then add a second vertex and join it to those of the 13 other vertices that have even degree.
(If the result is not a connected anti-Eulerian graph, it is rejected.)
By an argument similar to Theorem~\ref{thm:extend}, it can be shown that all graphs corresponding to
Type~II codes of length 14 must be part of this extended set.
Classifying all Type~II codes of length 14 turned out to be infeasible with our
computational resources. Even when using partitioning by partial weight distribution, the 
largest partitions were too large to be processed.
However, we were able to generate the LC orbits of all graphs corresponding to $(14,2^{14},6)$ codes.
Extremal Type~II codes have a unique weight enumerator, and the weight enumerator of a $(14,2^{14},6)$
code must be
\[
W_{14}(x,y) = x^{14} + 273 x^8 y^6 + 2457 x^6 y^8 + 7098 x^4 y^{10} + 6006 x^2 y^{12} + 549 y^{14}.
\]
Table~\ref{tab:14codes} lists the number of codes by automorphism group size. Note that
codes with 21, 168, and 2184 automorphisms were previously unknown.
Generator matrices of the codes are available on-line~\cite{database}.

\begin{table}
\centering
\caption{Number of $(14,2^{14},6)$ Codes with Automorphism Group of Size~$a$}
\label{tab:14codes}
\begin{tabular}{cc}
\toprule
$a$ &  \\
\midrule
1 & 625\\
2 & 258\\
3 &  27\\
4 &  38\\
6 &  27\\
8 &  13\\
12 &  7\\
18 &  1\\
21 &  1\\
24 & 16\\ 
28 &  1\\
36 &  1\\
48 &  1\\ 
84 &  1\\
168 & 1\\
2184 & 1\\
6552 & 1\\
\midrule
All & 1020 \\
\bottomrule
\end{tabular}
\end{table}

As mentioned before, the set of self-dual \emph{linear} codes over $\GF(4)$ is
a subset of the self-dual additive codes of Type~II.
Note that conjugation of single coordinates does not preserve the linearity of a code.
It was shown by Van~den~Nest~\cite{nestthesis} that the code $\mathcal{C}$ generated
by a matrix of the form $\Gamma + \omega I$ can not be linear. However, if there is a linear code 
equivalent to $\mathcal{C}$, it can be found by conjugating some coordinates.
Conjugating coordinates of $\mathcal{C}$ is equivalent to setting some diagonal elements of $\Gamma$ to 1.
Let $A$ be a binary diagonal matrix such that $\Gamma + A + \omega I$ generates
a linear code. Van~den~Nest~\cite{nestthesis} proved that 
$\mathcal{C}$ is equivalent to a linear code if and only if there exists such a matrix $A$ that satisfies
$\Gamma^2 + A \Gamma + \Gamma A + \Gamma + I = 0$.
A similar result was found by Glynn~et~al.~\cite{glynnbook}.
Using this method, it is easy to check whether the LC orbit of a given graph corresponds to
a linear code. However, self-dual linear codes over $\GF(4)$ have already been classified up to length 16, and
we have not found a way to extend this result using the graph approach.

We remark that if $\mathcal{C}$ is a self-dual additive code over $\GF(4)$
with generator matrix $\Gamma + \omega I$,
it can be shown that the additive code over $\mathbb{Z}_4$ generated by $2\Gamma + I$ has
the same weight distribution as $\mathcal{C}$.
It has also been shown~\cite{rains} that self-dual additive codes over $\GF(4)$ can be 
mapped to \emph{isodual} binary linear codes, i.e., codes that are equivalent to their duals,
by the mapping $0 \mapsto 00$, $1 \mapsto 11$, $\omega \mapsto 01$ and $\omega^2 \mapsto 10$.
A code over $\mathbb{Z}_4$ and a binary code obtained from the same self-dual additive code over 
$\GF(4)$ by these two methods are related by the well-known \emph{Gray map}.
There are also severals mappings from self-dual additive codes over 
$\GF(4)$ to self-dual and self-orthogonal binary linear codes~\cite{hohn,gaborit,project}.

An interesting problem, posed by Höhn~\cite{hohn}, is to find the smallest 
code with trivial automorphism group, i.e., automorphism group of size 1. 
We find that there is no such code of 
length up to 8, but there is a single code of length 9 with trivial automorphism group.
This code has generator matrix
\[
\left(
\begin{array}{ccccccccc}
\omega&0&0&0&0&0&0&1&1\\
0&\omega&0&0&1&0&0&1&0\\
0&0&\omega&1&0&0&1&0&0\\
0&0&1&\omega&0&0&0&1&1\\
0&1&0&0&\omega&1&0&0&1\\
0&0&0&0&1&\omega&1&0&1\\
0&0&1&0&0&1&\omega&0&1\\
1&1&0&1&0&0&0&\omega&0\\
1&0&0&1&1&1&1&0&\omega
\end{array}
\right).
\]
The smallest Type~II codes with trivial automorphism groups have length 12.
One such code is generated by
\[
\left(
\begin{array}{cccccccccccc}
\omega&0&0&0&0&0&0&1&0&0&1&1\\
0&\omega&0&0&0&0&1&0&1&0&0&1\\
0&0&\omega&0&1&1&0&0&0&0&0&1\\
0&0&0&\omega&1&0&1&0&0&0&1&0\\
0&0&1&1&\omega&0&0&1&0&0&0&0\\
0&0&1&0&0&\omega&0&0&1&1&1&1\\
0&1&0&1&0&0&\omega&0&0&1&1&1\\
1&0&0&0&1&0&0&\omega&1&1&0&1\\
0&1&0&0&0&1&0&1&\omega&1&0&1\\
0&0&0&0&0&1&1&1&1&\omega&1&0\\
1&0&0&1&0&1&1&0&0&1&\omega&0\\
1&1&1&0&0&1&1&1&1&0&0&\omega
\end{array}
\right).
\]
Table~\ref{tab:trivial} lists the numbers of Type~I and Type~II codes with
trivial automorphism group by length and minimum distance.
Note that for length 12, almost half the codes have trivial automorphism group.
For high lengths, one can expect almost all codes to have trivial automorphism group~\cite{hohn}.
This implies that the bound on $t_n$ given by Eq.~\eqref{eq:bound} is tighter for
higher $n$.
Observe that in Table~\ref{tab:trivial}, no code of minimum distance less than 3 is listed. 
It is easy to show that all codes with minimum distance 1 or 2 must
have nontrivial automorphisms.

\begin{table}
\centering
\caption{Number of Type~I (Type~II) Codes of Length~$n$ and Minimum Distance~$d$
with Trivial Automorphism Group}
\label{tab:trivial}
\begin{tabular}{cr@{ }rr@{ }rr@{ }rr@{ }rr@{ }rr@{ }r}
\toprule
$d \backslash n$ & \multicolumn{2}{c}{$\le8$} & \multicolumn{2}{c}{9} & \multicolumn{2}{c}{10} 
 & \multicolumn{2}{c}{11} & \multicolumn{2}{c}{12} & \multicolumn{2}{c}{14}\\
\midrule
3     &   &     & 1 & (0) & 113 & (0) & 6,247 & (0) & 392,649 & (0)   &?&(0)\\
4     &   &     &   &     &   3 & (0) & 1,180 & (0) & 163,982 & (102) &?&(?)\\
5     &   &     &   &     &     &     &       &     & 25      & (0)   &?&(0)\\
6     &   &     &   &     &     &     &       &     &         &       &?&(625)\\
\midrule
All   & 0 & (0) & 1 & (0) & 116 & (0) & 7,427 & (0) & 556,656 & (102) &?&(?)\\
\bottomrule
\end{tabular}
\end{table}

\section{Conclusions}\label{sec:concl}

By using graph representation and equivalence via local complementation,
we have classified all additive codes over $\GF(4)$ of length up to 12 
that are self-dual with respect to the Hermitian trace inner product. 
It follows from the bound given by Eq.~\eqref{eq:bound} that there
are at least 72,573,549 codes of length 13. It is not feasible
to classify all codes of length 13 using our method and the
computational resources available to us.
We were however able to classify the 1,020 extremal Type~II $(14,2^{14},6)$ codes.
This was done by exploiting the fact that Type~II codes correspond to
anti-Eulerian graphs. Finally, we showed that the smallest Type~I and Type~II codes 
with trivial automorphism group have length 9 and 12, respectively.

The graph representation of a self-dual additive code over $\GF(4)$ can
also give information about the properties of the code.
Tonchev~\cite{tonchev} showed that \emph{strongly regular} graphs
give rise to interesting codes. In particular, codes represented by the 
strongly regular \emph{Paley graphs} are well-known \emph{quadratic residue codes}.
We have shown that many extremal and optimal 
codes can be represented by \emph{nested regular graphs}~\cite{mscthesis,setapaper}.
Glynn~et~al.~\cite{glynnbook} showed that the minimum distance of a code
is equal to one plus the minimum \emph{vertex degree} over all
graphs in the corresponding LC orbit.
We have shown that the LC orbit corresponding to a code with high minimum distance only
contains graphs with both small \emph{independent sets} and small \emph{cliques}~\cite{mscthesis,setapaper}.

\paragraph*{Acknowledgements}
We would like to thank Philippe Gaborit for his helpful comments.
Also thanks to the Bergen Center for Computational Science, whose
cluster computer made the results in this paper possible.
This research was supported by the Research Council of Norway.

\end{document}